\documentclass[11pt,reqno]{amsproc}

\usepackage[margin=1in]{geometry}
\usepackage{mathrsfs}

\usepackage[usenames,dvipsnames]{color}
\usepackage[colorlinks=true, pdfstartview=FitV, linkcolor=RoyalBlue,
            citecolor=ForestGreen, urlcolor=blue]{hyperref}
\usepackage{fourier}

\title[Global Well-Posedness of a 3D MHD Model in Porous Media]
{Global Well-Posedness of a 3D MHD Model in Porous Media
\footnote{{\em Date:} June 4, 2018} }

\author[E. S. Titi]{Edriss S. Titi}
\address{\textnormal{(Edriss S. Titi)} Department of Mathematics\\
Texas A\&M University\\ College Station, TX 77843, USA.}
\email[E. S. Titi] {titi@math.tamu.edu}

\author{Saber Trabelsi}
\address{ Department of Mathematics \& Statistics, King Fahd University of Petroleum and Minerals, Dhahran 31261, KSA.}
\email{saber.trabelsi@kfupm.edu.sa}

\usepackage[pdftex]{graphicx}

\usepackage{amsmath, amsthm, amssymb,amsfonts,latexsym,verbatim,amsbsy}

\usepackage[usenames,dvipsnames]{color}
\usepackage{amsfonts}
\usepackage{mathrsfs}
\usepackage{ dsfont}
\usepackage{times}
\usepackage{hyperref}
\numberwithin{equation}{section}

\makeatletter
\@addtoreset{equation}{section}
\makeatother

\theoremstyle{plain}\newtheorem{thm}{Theorem}[section]
\theoremstyle{plain}
\theoremstyle{plain}
\theoremstyle{plain}
\theoremstyle{plain}\newtheorem{remark}[thm]{Remark}
\theoremstyle{plain}
\theoremstyle{plain}
\newcommand{\T}{{\Omega}}
\newcommand{\nr}{{|\!|}}
\newcommand{\iT}{{\int_{\Omega}\,}}
\newcommand{\dt}{{\,\frac{d}{dt}}\,}
\newcommand{\un}{{\,\nr u\nr}}

\newcommand{\ugn}{{\,\nr \nabla u\nr}}
\newcommand{\udn}{{\,\nr \Delta u\nr}}
\newcommand{\bn}{{\,\nr b\nr}}
\newcommand{\bgn}{{\,\nr \nabla b\nr}}
\newcommand{\bdn}{{\,\nr \Delta b\nr}}
\newcommand{\vn}{{\,\nr v\nr}}
\newcommand{\vgn}{{\,\nr \nabla v\nr}}

\begin{document}
\maketitle
\begin{abstract}
In this paper we show the global well-posedness of solutions to a three-dimensional  magnetohydrodynamical (MHD) model in porous media. Compared to the classical MHD equations, our system involves a nonlinear damping term in the momentum equations due to the ``Brinkman-Forcheimer-extended-Darcy'' law of flow in porous media.
\end{abstract}
\vskip6pt
\noindent \textbf{MSC class:} 76W05, 76S05, 35Q30, 35Q35, 76B03, 93C10, 93C20, 76B75.

\setcounter{tocdepth}{2}
\tableofcontents
\begin{center} \it This work is dedicated to Professor Darryl  Holm on the occasion of his 70$^{th}$ birthday  \end{center}
\section{Introduction}
\noindent The magnetohydrodynamics (MHD) equations form the system that governs  the interaction of electrically conducting fluids  and electromagnetic forces. They arise in various disciplines and applications ranging from the prediction of the space weather, exploratory geophysics, hydrology to the design of cooling system instruments and MHD generators. On the  modeling level, the MHD model involves a coupling, through the Lorentz force and Ohm's law for moving electrical conductors,  of the equations of fluid dynamics and the equations of electrodynamics. Briefly speaking, the MHD system consists of coupling  the Navier-Stokes and Maxwell's equations. Various physical situations require, sometimes, a modification or a simplification of these equations in order to capture the underlying physical phenomena at the relevant scales.
\vskip6pt\noindent
In this paper, we consider the motion of conducting fluid in porous medium. Nowadays, it is rather common to use the Darcy Law in the modeling of fluids momentum balance, i.e., the fluid motion, through the porous medium. Darcy's empirical law represents a simple linear relationship between the flow rate and the pressure drop in a porous medium, i.e.,
\begin{equation}\label{Darcy} u_f=- \frac k\mu\, \nabla p, \end{equation}
where $u_f$ is the Darcy velocity, $k$ is the permeability of the porous medium, $\mu$ is the dynamic viscosity of the fluid, and $p$ is the pressure. Any deviation from this scenario of flow in porous media is termed non-Darcy flow. Observe that Darcy's law neglects the inertia or the acceleration forces in the fluid when compared to the classical Navier-Stokes equations. Indeed, it is assumed that in a porous medium a large body with large surface area of the pores, is exposed to the fluid flow so that the viscous resistance will greatly exceed acceleration forces in the fluid unless turbulence sets in.  However, several physical situations  (e.g., in the case of relatively high velocity,  in the presence of important molecular and ionic effects, in the presence of non-Newtonian fluids etc.) violate the linearity relation of  Darcy's law (see, e.g., \cite{dupuit1863etudes,mcclure2010beyond,panfilov2006physical,fourar2004non} and references therein). An alternative approach is to use the Forchheimer correction of the Darcy law. More precisely, instead of \eqref{Darcy}, one may uses the Darcy-Forchheimer law \cite{forchheimer1901wasserbewegung},  given by
\begin{equation*} \nabla p= -\frac{\mu}{k}\,{\mathbf{ v}}_f- \gamma \rho_f |{\mathbf{ v}}_f|^2 \,{\mathbf {v}}_f,\end{equation*}
where $\gamma>0$ is the Forchheimer coefficient and ${\mathbf {v}}_f$ stands for the Forchheimer velocity and $\rho_f$ the density. In other words, Forchheimer law assumes that Darcy's law is still valid up to an additional nonlinear term  to account for the increased pressure drop for large velocity values. Then the Brinkman-Forcheheimer-extended-Darcy  equations (a generalisation actually) read
\begin{equation} {\label{BFeD}}
\left\lbrace\begin{array}{ll}
&\partial_t\,u-\nu\,\Delta\,u +(u\cdot\nabla)\,u +a\,|u|^{2\alpha}u+b\,|u|^{2\beta}u+\nabla\,p=f,\\ \\
& \nabla\cdot u=0,\,u\vert_{t=0}=u_0.
\end{array}
\right.
\end{equation}
This model was originally derived in its classical configuration ($\alpha=\frac12, \beta=0, a>0$, and $b>0$) in the framework of thermal dispersion in a porous medium using the method of volume averaging of the velocity and temperature deviations in the pores, see, e.g., \cite{HSU19901587}. A discussion of the formulation, validity and limitation of the BFD system can be found in \cite{nield1991limitations,vafai1987analysis}. The continuous dependence of BFeD equations on the Forchheimer coefficient with Dirichlet boundary conditions is studied in  \cite{louaked2015pseudocompressibility} (see also \cite{celebi2006continuous,liu2007structural,louaked2017approximation,payne1999convergence,straughan2008stability} and references therein, and \cite{celebi2005continuous} where a coupling with the temperature is considered). The long-time behavior of the solutions and the existence of global attractors to the BFeD system have been studied in \cite{ouyang2009note,uugurlu2008existence,wang2008existence,you2012existence} (see also \cite{cai2008weak} where the BFeD system is considered on $\mathbb R^3$ with the standard decay condition $|u| \to 0$, as $|x|\to \infty$) for restrictive range of the exponent $\alpha$. Also, existence, decay rates and some qualitative properties of weak solutions are shown in \cite{antontsev2010navier}.  In \cite {kalantarov2011smooth, markowich2016continuous}, the  BFeD system is investigated, and  existence and uniqueness of weak and strong solutions are shown. The authors of \cite{kalantarov2011smooth} assume Dirichlet boundary conditions and regular enough initial data. In \cite{markowich2016continuous}, the authors assume periodic boundary conditions (their result can be extended to the same boundary conditions as of \cite{kalantarov2011smooth} at the price of heavy technicalities due to the use of the maximal regularity estimates for the semi-linear stationary Stokes operator), but obtain the existence and uniqueness of strong solution with initial data less regular than in \cite{kalantarov2011smooth}. Eventually, an anisotropic viscous version of the BFeD system, and a 3D Forchheimer-B\'enard convection system are studied in \cite{bessaih2016existence} and  \cite{TrabMHD}, respectively.

\vskip6pt
\noindent
 In this paper, we consider the following three-dimensional MHD model of the motion of electrically conducting fluid in a porous medium
\begin{align*}
\mathcal {S} :\quad \left\lbrace\begin{array}{ll}
&\partial_t\,u-\nu\,\Delta\,u + u\cdot \nabla u- b\cdot \nabla b+a\,|u|^{2\alpha}u+\nabla p=0 ,\\ \\
&\partial_t\,b-\kappa\,\Delta\,b + u\cdot \nabla b- b\cdot \nabla u = 0,\\ \\
&\nabla\cdot u=0,\, u_{\vert_{t=0}}=u_0,\; b_{\vert{{t=0}}}=b_0,
\end{array}
\right.
\end{align*}
with periodic boundary conditions with fundamental periodic domain $\Omega=[0,L]^3$, $L>0$.  In system $\mathcal S$, the unknowns are $u$ the fluid velocity, $b$ the magnetic field, and $p$ the pressure. The viscosity  $\nu\geq 0$ and the magnetic resistivity $\kappa\geq0$, are given constants. The damping coefficients $a$ and the exponent $\alpha$ are nonnegative constants. The first equation $\mathcal S_1$ reflects the conservation of momentum of the fluid through the periodic porous medium with basic domain $\Omega$. On the physical level, this equation is believed to be accurate when the flow velocity is too large for Darcy's law to be valid. Equation $\mathcal S_2$ is the magnetic field evolution equation. In the sequel, we will consider divergence-free initial magnetic field $b_0$, that is $\nabla\cdot b_0=0$. Consequently, one concludes from  the second equation $\mathcal S_2$ that
\begin{align*}
\left\lbrace\begin{array}{ll}
&\partial_t (\nabla\cdot b) -\kappa\,\Delta \,(\nabla\cdot b) + (u\cdot\nabla)\, \nabla\cdot b=0,\\ \\
&(\nabla\cdot b)(0,x)=\nabla \cdot b_0(x).
\end{array}
\right.
\end{align*}
The above implies, for reasonably smooth $u$, that
\[\nabla\cdot b(t,x)=\nabla\cdot b_0(x)=0.\]
Hence, $\nabla\cdot b =0$ is a property of the solution, unlike the divergence-free velocity, $\nabla\cdot u =0$, which is a constraint.
\vskip6pt
\noindent
When $a=0$, system $\mathcal S$ reduces to the classical MHD equations. There is a very rich literature dedicated to the mathematical analysis of this system. In the two-dimensional viscous and resistive case (namely fully dissipative in fluids and in magnetic field, $\nu,\kappa > 0$),  the global existence of strong solutions is shown in \cite{Duvaut1972}. Some results about short-time existence and uniqueness of  strong solutions in the three-dimensional case are shown in \cite{Duvaut1972} and \cite{sermange1983some}. There are several interesting regularity criterion for system $\mathcal S$ (when $a=0$) in the three-dimensional case mainly based on Prodi-Serrin's condition. We refer the reader to, e.g., \cite{he2005regularity,he2005partial,cao2010two,kang2009interior} and references therein.
\vskip6pt\noindent
In the sequel, we shall intensively and implicitly use Young's inequality
\[rs\leq \epsilon r^p + \epsilon^{-q/p}\,s^q,\quad r,s>0,\quad 1/p+1/q=1,\quad \text{for}\quad p\geq 1, q\geq 1\quad \text{and any } \quad \epsilon>0.\]
Also, we will use  the Poincar\'e inequality
\begin{align*}
&\sqrt\lambda\,\nr \varphi\nr_{2} \leq \,\nr \nabla \varphi\nr_2, \quad \sqrt\lambda=\frac{2\pi}{L},
\end{align*}
 valid for every periodic function $\varphi \in H^1(\Omega)$  with zero mean in $\Omega$. Eventually, we shall use the following notation for the mean of a function $\varphi$
\[ \overline \varphi=\frac1{|\Omega|}\,\iT \varphi (x)\,dx.\]
Let us mention that we will use these notations indiscriminately for both scalars and vectors, which should not form any source of confusion.
\vskip6pt\noindent
Now, we are able to state our first result as follows
\begin{thm}\label{thm1}
Let $(u_0,b_0) \in L^2(\Omega)\times L^2(\Omega)$ such that $\nabla\cdot u_0=\nabla\cdot b_0=0$, and $\overline b_0=0$. Let $a,\nu,\kappa>0$ and $\alpha\geq0$. Then system $\mathcal S$ has  global weak solutions $(u(t),b(t))$ satisfying
\begin{align*}
&u(x,t) \in L^\infty(\mathbb R^+;L^2(\Omega))\cap L^2_{\rm loc}([0,+\infty); H^1(\Omega)) \cap L^{{2\alpha+2}}_{\rm loc}([0,+\infty);L^{{2\alpha+2}}(\Omega)),\\&{\rm and}\\
& b(x,t) \in L^\infty(\mathbb R^+;L^2(\Omega))\cap L^2_{\rm loc}([0,+\infty); H^1(\Omega)),\quad \text{with}\quad \nabla\cdot b(t)=0.
\end{align*}
In particular, the following holds true
\[ \limsup_{t\to+\infty}\,\left(\nr u(t)\nr_2^2 +\nr b(t)\nr_2^2\right)\leq  L^3\,\frac{a^{-\frac1{2\alpha-1}}}{\max\left\{\kappa,\nu\right\}}\,\left(\frac{\nu}{\lambda}\right)^{\frac{2\alpha+2}{2\alpha-1}}.\]
Moreover, if $\alpha\geq\frac32$ and $b_0 \in L^{3\frac{\alpha+1}{\alpha}}(\Omega)$, then the solution $(u(t),b(t))$ is unique, and depends continuously on the initial data.
\end{thm}
\vskip6pt \noindent Concerning strong solutions to system $\mathcal S$, we have the following:
\begin{thm}\label{thm2}
Let $(u_0,b_0) \in H^1(\Omega)\times H^1(\Omega)$ such that  $\nabla\cdot u_0=\nabla\cdot b_0=0$, and $\overline b_0=0$. Let $ a,\nu,\kappa>0$ and $\alpha\geq\frac32$. Then system $\mathcal S$ has  a unique global strong  solutions $(u(t),b(t))$ satisfying
\begin{align*}
&u(x,t) \in C_{\rm b}^0([0,+\infty);H^1(\Omega))\cap L^2_{\rm loc}([0,+\infty); H^2(\Omega)) \cap L^{{2\alpha+2}}_{\rm loc}([0,+\infty);L^{{2\alpha+2}}(\Omega)),\\&{\rm and}\\
& b(x,t) \in C_{\rm b}^0([0,+\infty);H^1(\Omega))\cap L^2_{\rm loc}([0,+\infty); H^2(\Omega)),\quad \text{with}\quad \nabla\cdot b(t)=0.
\end{align*}
Moreover, the solution depends continuously on the initial data.  Furthermore, if $u_0\in L^{2\alpha+2}(\Omega)$, then
\[\partial_t u(x,t) \in L^2_{\rm loc}([0,+\infty);L^2(\Omega)),\quad u\in L^{\infty}_{\rm loc}([0,+\infty);L^{{2\alpha+2}}(\Omega)), \quad \text{and}\quad\partial_t b(x,t) \in L^2_{\rm loc}([0,+\infty);L^2(\Omega)).\]
\end{thm}
\section{Global well-posedness}
\noindent  In this section, we prove Theorem \ref{thm1} and Theorem \ref{thm2}.
The global (in time) existence of solutions is shown in a classical way by proceeding in three steps. First, we use a Faedo-Galerkin approximation procedure to show the short time existence of solutions. The reader is referred, e.g., to any textbook about Navier-Stokes equations (e.g., \cite{constantin1988navier,temam2012infinite}) for details. Next, we obtain the necessary {\it a priori} bounds that allow to extend the solution of the Faedo-Galerkin system  globally in time. Eventually, we pass to the limit in the approximation procedure using the Aubin compactness Theorem, relying on the established {\it a priori} bounds for the solution and their derivatives (see, e.g., the details concerning similar systems in \cite{markowich2016continuous}).
\subsection{Existence of solutions} \label{existence sec}
The short-time existence of approximate solutions can be obtained through the standard Faedo-Galerkin approximation method, see, for instance, \cite{constantin1988navier,temam2012infinite} (see also \cite{markowich2016continuous}). Thus, we will omit this part and focus on the formal estimates.
\subsection{{ A Priori}  Estimates} \label{apriorisec}
In this section, we perform formal calculation on system $\mathcal{S}$ to obtain the needed {\it a priori} estimates.As we mentioned above, those calculations can be rigorously justified by using the Faedo-Galerkin approximation. In the sequel, we will assume that $\nu,\kappa, a>0$.
\subsubsection{$L^2(\Omega)$ estimates of the velocity and the magnetic field}
Since $b_0$  satisfies $\nabla \cdot b_0=0$, it follows that $\nabla\cdot b(t)=0$ for all $t\geq0$. Moreover,  from  equation $\mathcal S_2$, we get
\[\frac d{dt}\,\int_{\Omega}\,b(t,x)\,dx=0, \]
and therefore
\begin{equation}
\overline{b}(t)=\overline{b}_0=0,\quad \text{for all }\quad t\geq 0.
\label{propaverage}\end{equation}
On the other hand, from equation $\mathcal S_1$, we have
\[\frac d{dt}\,\int_{\Omega}\,u(t,x)\,dx= -a\int \,|u(t,x)|^{2\alpha}\,u(t,x) \,dx.\]
 Hence, in general it is not true that  $\overline {u}(t,x)=\overline {u}_0$.
 \vskip6pt\noindent Now, let $u_0, b_0 \in L^2(\Omega)$ such that  $\nabla\cdot u_0=\nabla\cdot b_0=0, \overline b_0=0$ and $\alpha\geq 0$. In this section we  will establish  $L^2(\Omega)$ bounds for $u(t)$ and $b(t)$. First, we multiply $\mathcal S_1$ by $u$ and integrate over $\T$ to obtain
\begin{align*}
\frac12\dt\un_2^2 +\nu\ugn_2^2+ a\un^{2(\alpha+1)}_{2(\alpha+1)} = \iT  b\cdot \nabla b\cdot u\,dx.
\end{align*}
Next, we multiply $\mathcal S_2$ by $b$ and integrate over $\T$, and integrate by parts to get
\begin{align*}
\frac12\dt\bn_2^2 +\kappa\bgn_2^2 = \iT  b\cdot \nabla u\cdot b\,dx=-\iT  b\cdot \nabla b\cdot u\,dx
\end{align*}
Summing-up the latter inequalities, we can write
\begin{align}\label{allbacktothis}
\dt \left(\nr u(t)\nr_2^2+\nr b(t)\nr_2^2\right)+2\min\{\nu,\kappa\} \vgn_2^2  + 2a\un^{2(\alpha+1)}_{2(\alpha+1)}=0.
\end{align}
Integrating this inequality with respect to time on $(0,t)$ leads to
\begin{align}\label{fewbounds}
\nr u(t)\nr_2^2+\nr b(t)\nr_2^2\leq \nr u_0\nr_2^2+\nr b_0\nr_2^2,\quad \int_0^{t}\,\left(\nr \nabla u(\tau)\nr_2^2 +\nr \nabla b(\tau)\nr_2^2 \right)\,d\tau\leq \frac{\nr u_0\nr_2^2+\nr b_0\nr_2^2}{2\min\{\nu,\kappa\}},
\end{align}
and
\begin{align*}
 \int_0^{t}\,\nr u(\tau)\nr^{2(\alpha+1)}_{2(\alpha+1)} \,d\tau\leq \frac{\nr u_0\nr_2^2+\nr b_0\nr_2^2}{2a}.
\end{align*}
Consequently, we have
\begin{equation}\label{propweaksolu}u(t)\in L^\infty(\mathbb R^+;L^{2}(\Omega))\,\cap\, L^2_{\rm loc}([0,+\infty);H^1(\Omega)) \,\cap\,L^{2\alpha+2}_{\rm loc}([0,+\infty);L^{2\alpha+2}(\Omega)), \end{equation}
and
\[b(t)\in L^\infty(\mathbb R^+;L^{2}(\Omega))\,\cap\, L^2_{\rm loc}([0,+\infty);H^1(\Omega)),\]
for all $\alpha\geq0$.
This bound can be slightly improved. Indeed, observe that using H\"older and Young inequalities, one has the following
\begin{equation}\label{interpolate2alphaplus2}
\frac{\nu}{\lambda}\,\un_2^2 \leq a\,\un_{{2\alpha+2}}^{2\alpha+2} + L^3\,a^{-\frac1{2\alpha-1}}\,\left(\frac{\nu}{\lambda}\right)^{\frac{2\alpha+2}{2\alpha-1}}.
\end{equation}
Next, thanks to Poincar\'e inequality and \eqref{allbacktothis}, we have
\begin{align*}
\dt\un_2^2+\dt\bn_2^2  +2\nu\ugn_2^2 +\frac{2\kappa}{\lambda}\bn_2^2 + \frac{2\nu}{\lambda}\,\un^2_2 \leq 2L^3\,a^{-\frac1{2\alpha-1}}\,\left(\frac{\nu}{\lambda}\right)^{\frac{2\alpha+2}{2\alpha-1}}.
\end{align*}
Thus, we have
\begin{align*}
\dt \left(\nr u(t)\nr_2^2+\nr b(t)\nr_2^2\right) + \frac2\lambda\,\max\left\{{\kappa},{\nu}\right\} \vn_2^2\leq 2L^3\,a^{-\frac1{2\alpha-1}}\,\left(\frac{\nu}{\lambda}\right)^{\frac{2\alpha+2}{2\alpha-1}},
\end{align*}
from which we infer that for all $t\geq0$
\begin{align*}\label{limsupv2}
\nr u(t)\nr_2^2+\nr b(t)\nr_2^2&\leq \left(\nr u_0\nr_2^2+\nr b_0\nr_2^2\right)\,e^{- \frac2\lambda\max\left\{{\kappa},\nu\right\}\,t} \nonumber\\
&+  L^3\,\frac{a^{-\frac1{2\alpha-1}}}{\max\left\{\kappa,\nu\right\}}\,\left(\frac{\nu}{\lambda}\right)^{\frac{2\alpha+2}{2\alpha-1}} \,\left(1-e^{- 2\max\left\{\frac{\kappa}{\lambda},\frac\nu\lambda\right\}\,t}\right).
\end{align*}
Consequently, we have
\[ \limsup_{t\to+\infty}\,\left(\nr u(t)\nr_2^2+\nr b(t)\nr_2^2\right)\leq   L^3\,\frac{a^{-\frac1{2\alpha-1}}}{\max\left\{\kappa,\nu\right\}}\,\left(\frac{\nu}{\lambda}\right)^{\frac{2\alpha+2}{2\alpha-1}}.\]
In particular, system $\mathcal S$ has an absorbing ball in $L^2(\Omega)\times L^2(\Omega)$, provided $\alpha\geq 0$.
\subsubsection{$L^{3\frac{\alpha+1}{\alpha}}(\Omega)$ estimate of the magnetic field}\label{mainsec}
In this section we establish {\it a priori} bound for the magnetic field in the $L^{3\frac{\alpha+1}{\alpha}}(\Omega)$ norm, provided that the initial data $u_0, b_0 \in L^2(\Omega)$ such that  $\nabla\cdot u_0=\nabla\cdot b_0=0,  \overline b_0=0$, and in addition $b_0\in L^{3\frac{\alpha+1}{\alpha}}(\Omega)$ and $\alpha\geq\frac32$.
\vskip6pt\noindent
  We multiply $\mathcal S_2$ by $ |b|^{\frac{\alpha+3}{\alpha}}b$ and integrate over $\T$ to obtain
\begin{align*}
\nonumber\frac13\frac\alpha{\alpha+1}\dt \bn_{3\frac{\alpha+1}{\alpha}}^{3\frac{\alpha+1}{\alpha}} - \kappa\,\iT\Delta\,b\cdot b\,|b|^\frac{\alpha+3}{\alpha}\,dx &=  \iT  (b\cdot \nabla) \,u\cdot  b\,|b|^\frac{\alpha+3}{\alpha}\,dx\\
& =-\iT  b\cdot \nabla(|b|^\frac{\alpha+3}{\alpha} b) \cdot u \,dx.
\end{align*}
On the one hand, by using  Stroock-Varopoulos inequality (see, e.g., \cite{liskevich96}), we have
\begin{align}\label{Varpolous}
- \kappa\,\iT\Delta\,b\cdot b\,|b|^\frac{\alpha+3}{\alpha}\,dx \geq \frac{4\kappa}{9}\,\frac{\alpha(2\alpha+3)}{(\alpha+1)^2}\,\nr \nabla|b|^{\frac{3}{2}\frac{\alpha+1}{\alpha}}\nr_2^2 \geq \frac{5\kappa}{9}\,\nr \nabla|b|^{\frac{3}{2}\frac{\alpha+1}{\alpha}}\nr_2^2.
\end{align}
On the other hand, it is rather easy to see that
\begin{align*}
\left|\iT  b\cdot \nabla(|b|^\frac{\alpha+3}{\alpha} b) \cdot u \,dx\right|&\leq \frac{\alpha+3}{3(\alpha+1)}\,\iT\,|b|^{\frac32\frac{\alpha+1}{\alpha}}\,|\nabla|b|^{\frac32\frac{\alpha+1}{\alpha}}|\,|u|\,dx \\
&\leq \frac23\,\iT\,|b|^{\frac32\frac{\alpha+1}{\alpha}}\,|\nabla|b|^{\frac32\frac{\alpha+1}{\alpha}}|\,|u|\,dx \\
\nonumber\text{\tiny{[By Cauchy-Schwarz and Young  ineq.]}} &\leq \frac{2\kappa}9\,\nr \nabla|b|^{\frac{3}{2}\frac{\alpha+1}{\alpha}}\nr_2^2+\frac 2\kappa \,\nr u|b|^{\frac{3}{2}\frac{\alpha+1}{\alpha}}\nr_2^2.
\end{align*}
Now, using H\"older and Young  inequalities, we obtain
\begin{align*}
\nr u|b|^{\frac{3}{2}\frac{\alpha+1}{\alpha}}\nr_2^2 &\leq \nr u\nr^2_{2\alpha+2}\,\nr  |b|^{\frac{3}{2}\frac{\alpha+1}{\alpha}}\nr_{2\frac{\alpha+1}{\alpha}}^{2} \\
\nonumber\text{\tiny{[By Gagliardo-Nirenberg  ineq.]}} &\leq C\, \nr u\nr^2_{2\alpha+2}\,\nr  |\nabla |b||^{\frac{3}{2}\frac{\alpha+1}{\alpha}}\nr_{2}^{\frac{3}{\alpha+1}}\,\nr b\nr^{3-\frac3{2\alpha}}_{3\frac{\alpha+1}{\alpha}} \\
\nonumber\text{\tiny{[By Young's  ineq.]}} &\leq \frac {\kappa^2}9 \, \nr \nabla|b|^{\frac{3}{2}\frac{\alpha+1}{\alpha}}\nr_2^2 + C\,\kappa^{\frac6{1-2\alpha}}\,\nr u\nr_{2\alpha+2}^{\frac{4\alpha+4}{2\alpha-1}}\, \nr b\nr_{3\frac{\alpha+1}{\alpha}}^{3\frac{\alpha+1}{\alpha}},
\end{align*}
where $C$ denotes a positive constant dependent on $L$. Now, if we assume $\alpha\geq \frac32$, then we can write
\[\nr u\nr_{2\alpha+2}^{\frac{4\alpha+4}{2\alpha-1}}\leq 1+\nr u\nr_{2\alpha+2}^{2\alpha+2}.\]
All in all, we obtain for all $\alpha\geq \frac32$
\begin{align*}
\dt \bn_{3\frac{\alpha+1}{\alpha}}^{3\frac{\alpha+1}{\alpha}} + \,\frac{\kappa}{3} \,\nr \nabla|b|^{\frac{3}{2}\frac{\alpha+1}{\alpha}}\nr_2^2 \leq  C\,\kappa^{\frac6{1-2\alpha}}\,(1+\nr u\nr_{2\alpha+2}^{2\alpha+2})\, \nr b\nr_{3\frac{\alpha+1}{\alpha}}^{3\frac{\alpha+1}{\alpha}}.
\end{align*}
Thanks to \eqref{fewbounds}, we  have $u(t,x) \in L^{2\alpha+2}_{\rm loc}([0,+\infty); L^{2\alpha+2}(\Omega))$, thus
\begin{align*}
\nr b(t)\nr_{3\frac{\alpha+1}{\alpha}} &\leq \nr b_0\nr_{3\frac{\alpha+1}{\alpha}}\exp\left\{ C\,\kappa^{\frac6{1-2\alpha}}\,\int_0^t\,\left(1+ \nr u(\tau)\nr^{2\alpha+2}_{2\alpha+2}\right)\,d\tau \right\} \\
&\leq \nr b_0\nr_{3\frac{\alpha+1}{\alpha}}\exp\left\{ C\,\kappa^{\frac6{1-2\alpha}}\, \left(t+\frac{\nr u_0\nr_2^2+\nr b_0\nr_2^2}{2a} \right)\right\}.
\end{align*}
In particular,
\begin{equation}
b(t)\in L^\infty_{\rm loc}([0,+\infty);L^{3\frac{\alpha+1}{\alpha}}(\Omega)),
\label{bregfortimeder}\end{equation}
provided $\alpha\geq \frac32$.
\subsubsection{$H^1(\Omega)$ estimates of the velocity and the magnetic field}
In this section, we establish  {\it a priori} estimates for  the velocity and the magnetic fields in $H^1(\Omega)$, provided that  $u_0, b_0 \in H^1(\Omega)$ with  $\nabla\cdot u_0=\nabla\cdot b_0=0, \overline b_0=0$, and $\alpha \geq \frac32$.

\vskip6pt
\noindent We multiply $\mathcal S_1$ by $-\Delta u$ and integrate over $\T$ to obtain
\begin{align}\label{strongtest}
\frac12\dt\ugn_2^2 +\nu\udn_2^2- a\,\iT |u|^{2\alpha}u\cdot \Delta u=  \iT  (u\cdot \nabla)\, u\cdot \Delta u\,dx-\iT  (b\cdot \nabla)\, b\cdot \Delta u\,dx.
\end{align}
On the one hand,
\begin{equation}\label{powertermdirect}- a\,\iT |u|^{2\alpha}u\cdot \Delta u\,dx=a\,(1+2\alpha) \,\nr |u|^\alpha\,\nabla u\nr^2_2. \end{equation}
On the other hand, since $\alpha>1$ and use H\"older and Young inequalities  to get
\begin{align}\label{cumbersome}
\left|\iT  (u\cdot \nabla)\, u\cdot \Delta u\,dx \right| &\leq \iT\,|u|\,|\nabla\,u|^\frac1\alpha\,|\nabla\,u|^{1-\frac1\alpha}\,|\Delta\,u|\,dx \nonumber\\
&\leq \nr u\,|\nabla\,u|^\frac1\alpha\nr_{{2\alpha}}\,\nr |\nabla\,u|^{1-\frac1\alpha}\nr_{\frac{2\alpha}{\alpha-1}}\,\udn_2\nonumber\\
&\leq \frac 2{\nu}\, \nr|u|^{\alpha}\,\nabla\,u\nr^{\frac2\alpha}_{2}\,\ugn^{2(1-\frac1\alpha)}_2  +\frac\nu8\,\udn_2^2\nonumber\\
&\leq \,\frac{2\epsilon}{\nu\,}\, \nr|u|^{\alpha}\,\nabla u\nr^{2}_{2} + \,\frac{2\epsilon^{\frac1{1-\alpha}}}{\nu} \,\ugn_2^2  +\frac\nu8\,\udn_2^2,\quad \forall \epsilon>0.
\end{align}
Next, thanks to Cauchy-Schwarz and Young inequalities, we have
\begin{align}
\nonumber\left|\iT  (b\cdot \nabla)\, b\cdot \Delta u\,dx\right| &\leq \frac 2\nu \,\nr b |\nabla b|\nr_2^2+\frac\nu8\udn_2^2.
\end{align}
Hereafter, $c$ will denote a dimensionless positive constant that might vary  form line to line. By virtue of the H\"older inequality we have
\begin{align}\label{topintlater}
\nonumber\frac 2\nu\,\nr b|\nabla b|\nr_2^2 &\leq \frac 2\nu\,\nr b\nr^2_{3\frac{\alpha+1}\alpha}\,\nr \nabla b\nr^2_{6\frac{\alpha+1}{\alpha+3}} \\
\nonumber\text{\tiny{[By Gagliardo-Nirenberg  ineq.]}}&\leq \frac {c}\nu\,\nr b\nr^2_{3\frac{\alpha+1}\alpha}\,\nr \nabla b\nr^{\frac{2}{\alpha+1}}_{2}\,\nr \Delta b\nr_2^{\frac{2\alpha}{\alpha+1}} \\
\text{\tiny{[By Young's  ineq.]}}&\leq  \frac{c}{\nu^{\alpha+1}\,\kappa^\alpha}\,\nr \nabla b\nr^{2}_{2} + \frac\kappa {6}\,\nr \Delta b\nr_2^{2},
\end{align}
where we used, above, the following version of the Gagliardo-Nirenberg inequality for  any periodic function $\varphi\in H^2(\Omega)$,
\begin{equation}\nr \nabla \varphi\nr_{6\frac{\alpha+1}{\alpha+3}} \leq c\,\nr \nabla \varphi\nr^{\frac{1}{\alpha+1}}_{2}\,\nr \Delta \varphi\nr_2^{\frac{\alpha}{\alpha+1}} \label{GN}\end{equation}
Combining the latter two inequalities, we get
\begin{align}\label{torefinelater}
\left|\iT  (b\cdot \nabla)\, b\cdot \Delta u\,dx\right| &\leq  \frac{c}{\nu^{\alpha+1}\,\kappa^\alpha}\,\nr \nabla b\nr^{2}_{2} + \frac\kappa  {6}\,\nr \Delta b\nr_2^{2} +\frac\nu8\udn_2^2.
\end{align}
Now, we multiply $\mathcal S_2$ by $-\Delta b$ and integrate over $\T$
\begin{align}\label{toreuseannihilation}
\frac12\dt\bgn_2^2 +\kappa\bdn_2^2 = \iT  u\cdot \nabla b\cdot \Delta b\,dx-\iT  b\cdot \nabla u\cdot \Delta b\,dx.
\end{align}
We start by estimating the last term in the right-hand hand of the latter equation. On the one side, using  Cauchy-Schwarz and Young inequalities, we have
\begin{align}\label{thistermasthefirst}
\left|\iT  (b\cdot \nabla)\, u\cdot \Delta b\,dx\right| &\leq \frac 3{2\kappa} \,\nr b |\nabla u|\nr_2^2+\frac\kappa {6}\bdn_2^2.
\end{align}
On the other hand, we have
\begin{align}
\nonumber\frac 3{2\kappa}\,\nr b|\nabla u|\nr_2^2 &\leq \frac 3{2\kappa}\,\nr b\nr^2_{3\frac{\alpha+1}\alpha}\,\nr \nabla u\nr^2_{6\frac{\alpha+1}{\alpha+3}} \\
\nonumber\text{\tiny{[By Gagliardo-Nirenberg  ineq.]}}&\leq \frac {c}\kappa\,\nr b\nr^2_{3\frac{\alpha+1}\alpha}\,\nr \nabla u\nr^{\frac{2}{\alpha+1}}_{2}\,\nr \Delta u\nr_2^{\frac{2\alpha}{\alpha+1}} b\nr^2_{3\frac{\alpha+1}\alpha}\\
\text{\tiny{[By Young's  ineq.]}}&\leq  \frac {c}{\nu^{\alpha}\kappa^{\alpha+1}} \,\nr b\nr^{2\alpha+2}_{3\frac{\alpha+1}\alpha} \,\nr \nabla u\nr^{2}_{2} + \frac\nu 8\,\nr \Delta u\nr_2^{2},\label{tousefortimeinteg}
\end{align}
where we used the Gagliardo-Nirenberg inequality \eqref{GN} for the velocity. Thus, we infer the following estimate
\begin{align}\label{lastforgrad}
\left|\iT  (b\cdot \nabla)\, u\cdot \Delta b\,dx\right| &\leq     \frac {c}{\nu^{\alpha}\kappa^{\alpha+1}} \,\nr b\nr^{2\alpha+2}_{3\frac{\alpha+1}\alpha}\,\nr \nabla u\nr^{2}_{2} + \frac\nu 8\,\nr \Delta u\nr_2^{2}+ \frac\kappa {6}\bdn_2^2.
\end{align}
Now,  we handle the remaining term using integration by parts
\begin{align*}
\iT  (u\cdot \nabla)\, b\cdot \Delta b\,dx&= \sum_{k,l,m=1}^3\,\iT u_k\,\partial_k b_l\,\partial_m^2\,b_l\,dx \\
&= - \sum_{k,l,m=1}^3\,\iT \partial_m\,u_k\,\partial_k b_l\,\partial_m\,b_l\,dx -  \sum_{k,l,m=1}^3\,\iT u_k\,\partial_m\partial_k b_l\,\partial_m\,b_l\,dx \\
&= - \sum_{k,l,m=1}^3\,\iT \partial_m\,u_k\,\partial_k b_l\,\partial_m\,b_l\,dx - \sum_{m=1}^3\,\iT u\cdot \nabla \partial_mb\cdot \partial_m b\,dx.
\end{align*}
Clearly, the last term of the right-hand side of this inequality vanishes. Next, integrating by parts the first term and using the divergence free condition, we obtain
\begin{align*}
 - \sum_{k,l,m=1}^3\,\iT \partial_m\,u_k\,\partial_k b_l\,\partial_m\,b_l\,dx &= \sum_{k,l,m=1}^3\,\iT \partial^2_m\,u_k\,\partial_k b_l\,\,b_l\,dx  + \sum_{k,l,m=1}^3\,\iT \partial_m\,u_k\,\partial_m\,\partial_k b_l\,b_l\,dx \\
 &= \iT (\Delta u\cdot \nabla)\, b\cdot  b\,dx+ \sum_{k,l,m=1}^3\,\iT \partial_m\,u_k\,\partial_m\,\partial_k b_l\,b_l\,dx. \\
 &= \sum_{m=1}^3\,\iT\,(\partial_mu\cdot\nabla)\partial_mb\cdot b\,dx.
\end{align*}
Again the  the first term in the second line above vanishes. Eventually, from the above we have
\[\left|\iT  (u\cdot \nabla)\, b\cdot \Delta b\,dx \right|\leq \left|\iT  |\nabla u|\,|\Delta b|\,|b|\, dx\right|,\]
and therefore it can be estimated exactly as in \eqref{thistermasthefirst}. Thus, we get
\begin{align}\label{ouff}
\left|\iT  (u\cdot \nabla)\, b\cdot \Delta b\,dx \right|\leq    \frac {c}{\nu^{\alpha}\kappa^{\alpha+1}} \,\nr b\nr^{2\alpha+2}_{3\frac{\alpha+1}\alpha}\,\nr \nabla u\nr^{2}_{2} + \frac\nu 8\,\nr \Delta u\nr_2^{2}+ \frac\kappa {6}\bdn_2^2.
\end{align}
Now, gathering (\ref{strongtest}-\ref{ouff}), and setting $\epsilon= \frac{a(1+2\alpha)\nu}{2}$ in \eqref{cumbersome}, we obtain
\begin{align}\label{togetbacktoitinrem}
&\dt\ugn_2^2 +\dt\bgn_2^2+\nu\udn_2^2+\kappa\bdn_2^2 +a(1+2\alpha)\ \nr |u|^\alpha\nabla u\nr_2^2 \\
\nonumber&\leq c\,\left(\,(a(1+2\alpha)\,\nu^\alpha)^{\frac1{1-\alpha}} + \frac {1}{\nu^{\alpha}\kappa^{\alpha+1}} \,\nr b\nr^{2\alpha+2}_{3\frac{\alpha+1}\alpha} \right)\,\nr \nabla u\nr_2^2 +\frac{c}{\nu^{\alpha+1}\,\kappa^\alpha}\,\nr \nabla b\nr_2^2.
\end{align}
Since $H^1(\Omega) \hookrightarrow L^{3\frac{\alpha+1}\alpha}(\Omega)$ for all $\alpha\geq 1$, we have $b_0\in L^{3\frac{\alpha+1}\alpha}(\Omega)$. Now,  recalling that $b(t)\in L^\infty_{\rm loc}([0,+\infty);L^{3\frac{\alpha+1}{\alpha}}(\Omega))$, for all $\alpha\geq \frac32$ thanks to \eqref{bregfortimeder}, then Gronwall's inequality leads clearly to the fact that
\begin{align}\label{regularityub}&\nabla u,\nabla b \in L^\infty_{\rm loc}([0,+\infty), L^2(\Omega)), \quad \Delta u,\Delta b \in L^2_{\rm loc}([0,+\infty)), L^2(\Omega)), \\ &|u|^\alpha\nabla u \in L^2_{\rm loc}([0,+\infty), L^2(\Omega)), \nonumber
\end{align}
\vskip6pt
\noindent Now, we  show that under the same assumptions on the initial data and $\alpha$, we have
\begin{equation}\partial_t u,\,\partial_t b\in L^2_{\rm loc}([0,+\infty), L^2(\Omega)).\label{dtestimate}\end{equation}
For this purpose, we multiply $\mathcal S_1$ by $\partial_t\,u$ and integrate over $\T$ to obtain
\[\nr \partial_t u\nr^2_2 + \frac\nu2\,\frac d{dt}\,\nr \nabla u\nr_2^2 + \iT  (u\cdot \nabla)\, u\cdot  \partial_t u\,dx  + \iT  (b\cdot \nabla)\, b\cdot \partial_t u\,dx+\frac{2}{2\alpha+2}\,\frac d{dt}\,\nr u\nr_{2\alpha+2}^{2\alpha+2}\leq 0 . \]
Now, on the one hand,
\[  \left|\iT  (u\cdot \nabla)\, u\cdot  \partial_t u\,dx  \right| \leq   \nr u|\nabla u|\nr_2^2  +\frac14\,\nr \partial_t u\nr_2^2.\]
On the other hand, we have
\begin{align*}  \left|\iT  (b\cdot \nabla)\, b\cdot  \partial_t u\,dx  \right| &\leq   \nr b|\nabla b|\nr_2^2  +\frac14\,\nr \partial_t u\nr_2^2.
\end{align*}
Gathering these estimates together, we get
\begin{align*}
\nr \partial_t u\nr^2_2 + \nu\,\frac d{dt}\,\nr \nabla u\nr_2^2+ \frac{1}{2\alpha+2}\,\frac d{dt}\,\nr u\nr_{2\alpha+2}^{2\alpha+2} \leq   2\,\nr u|\nabla u|\nr_2^2  + 2\,\nr b |\nabla b|\nr_2^2.
\end{align*}
Integrating this inequality with respect to time, we obtain
\begin{align}\label{temptoshowdt}
\nonumber \int_0^t\,\nr \partial_t u(\tau)\nr^2_2\,d\tau + \nu\,\nr \nabla u(t)\nr_2^2 &+  \frac{1}{2\alpha+2}\,\nr u(t)\nr_{2\alpha+2} \leq  \nu\,\nr \nabla u_0\nr_2^2 +  \frac{1}{2\alpha+2}\,\nr u_0\nr_{2\alpha+2}^{2\alpha+2} \\
&+ 2\,\int_0^t\,\nr u(\tau)|\nabla u(\tau)|\nr_2^2\,d\tau  + 2\,\int_0^t\,\nr b(\tau)|\nabla b(\tau)|\nr_2^2\,d\tau.
\end{align}
First, using Young's inequality, it is rather easy to see that
\begin{align*} \nr u|\nabla u|\nr_2^2&\leq \int_\Omega\, |u|^2\,|\nabla u|^\frac2\alpha\,|\nabla u|^{2-\frac2\alpha}\,dx \leq \nu^{1-\alpha}\,\nr |u|^\alpha |\nabla u|\nr_2^2+\nu\,\nr \nabla u\nr_2^2.\label{part1forcompact}
\end{align*}
Next, the term $\nr |b|\nabla b\nr_2^2 $ can be  estimated as in \eqref{topintlater} (recalling \eqref{bregfortimeder}) or only using Agmon's inequality (see \cite{Agmon})
\[\nr b\nr_\infty \leq c \,\nr \Delta b\nr^{\frac34}_2\,\nr b\nr^{ \frac14}_2,\]
where $c$ denotes a nonnegative constant. Indeed, we have
\begin{equation}\nr b|\nabla b|\nr_2^2 \leq \nr \nabla b\nr_{2}^2 \,\nr b\nr_\infty^2\leq c\,\nr \nabla b\nr_{2}^2 \,\nr \Delta b\nr^{\frac32}_2\,\nr b\nr^{ \frac12}_2 \leq c\,\nr \nabla b\nr_{2}^8\,\nr b\nr_2^2 +\,\nr \Delta b\nr_{2}^2,\label{part2forcompact}\end{equation}
Therefore, thanks to \eqref{regularityub}, we obtain  that \[u|\nabla u|, b|\nabla b| \in L^2_{\rm loc}([0,+\infty), L^2(\Omega)).\]
Eventually, thanks to \eqref{temptoshowdt}, we infer the first part of \eqref{dtestimate}. In order to prove the second part, we multiply $\mathcal S_2$ by $\partial_t\,b$ and integrate over $\T$ to get
\begin{equation}\label{tointegfordtb}\nr \partial_t b\nr^2_2 + \frac\kappa2\,\frac d{dt}\,\nr \nabla b\nr_2^2 + \iT  (u\cdot \nabla)\, b\cdot  \partial_t b\,dx  - \iT  (b\cdot \nabla)\, u\cdot \partial_t b\,dx= 0. \end{equation}
Now,  using H\"older and Young inequalities, we have
\[\left|  \iT  (u\cdot \nabla b)\,\cdot  \partial_t b\,dx\right| \leq \nr u|\nabla b|\nr_2^2+\frac14\nr\partial_tb\nr_2^2.\]
and
\[\left|  \iT  (b\cdot \nabla)\, u\cdot  \partial_t b\,dx\right| \leq \nr b|\nabla u|\nr_2^2+\frac14\nr\partial_tb\nr_2^2.\]
Thus, integrating \eqref{tointegfordtb} with respect to time, we get
\begin{align}\label{integratedfordt}
\int_0^t\,\nr \partial_t b(\tau)\nr^2_2\,d\tau &+ \kappa\, \nr \nabla b(t)\nr_2^2  \leq  \kappa\,\nr \nabla b_0\nr_2^2 + 2\int_0^t\,\left( \nr u(\tau)|\nabla b(\tau)|\nr_2^2+\nr b(\tau)|\nabla u(\tau)|\nr_2^2\right)\,d\tau.
\end{align}
Now, using Gagliardo-Nirenberg inequality, we have using
\begin{align*}
\nr u|\nabla b|\nr_2^2 \leq \nr u\nr_\infty^2\,\nr \nabla b\nr_2^2&\leq c \,\nr \Delta u\nr_2^2\,\nr \nabla u\nr_2^2\,\nr \nabla b\nr_2^2 + c\,\nr u\nr_2^2\,\nr \nabla b\nr_2^2 \nonumber\\
\text{\tiny{[By Young's  ineq.]}}& \leq \nr \Delta u\nr_2^2+ c\,\nr \nabla u\nr_2^4\,\nr \nabla b\nr_2^4 + c\,\nr u\nr_2^2\,\nr \nabla b\nr_2^2. \label{fordtinteg1}
\end{align*}
Furthermore, using \eqref{tousefortimeinteg}, we have
\begin{align*}
 \nr b|\nabla u|\nr_2^2
&\leq  \frac {c}{\nu^{\alpha}\kappa^{\alpha}} \,\nr b\nr^{2\alpha+2}_{3\frac{\alpha+1}\alpha} \,\nr \nabla u\nr^{2}_{2} + \frac{\nu\kappa} {16}\,\nr \Delta u\nr_2^{2}.
\end{align*}
Thus, thanks to \eqref{bregfortimeder} and \eqref{regularityub}, we obtain
\[u|\nabla b|, b|\nabla u| \in L^2_{\rm loc}([0,+\infty), L^2(\Omega)).\]
All in all, \eqref{integratedfordt} provides  the second part of \eqref{dtestimate}.
\subsection{Continuous dependence on the initial data}
In this section, we show the continuous dependence of the weak and strong solutions on the initial data, in particular their uniqueness. Let $(u,p,b)$ and $(\tilde u,\tilde p,\tilde b)$  associated to the initial data $(u_0,p_0,b_0)$ and $(\tilde u_0,\tilde p_0,\tilde b_0)$, respectively satisfying the respected conditions in Theorems \ref{thm1} and \ref{thm2} . Therefore, if we let $(v,q,d)=(u-\tilde u,p-\tilde p,b-\tilde b)$, then  $(v,q,d)$  satisfies
\begin{align*}
\tilde{\mathcal S}:\:\quad \left\lbrace\begin{array}{ll}
&\partial_t\,v-\nu\,\Delta\,v + u\cdot \nabla v +v\cdot \nabla \tilde u- b\cdot \nabla d-d\cdot \nabla \tilde b+a\,|u|^{2\alpha}u-a\,|\tilde u|^{2\alpha}\tilde u+\nabla q=0 ,\\ \\
&\partial_t\,d-\kappa\,\Delta\,d + v\cdot \nabla b+ \tilde u\cdot \nabla d- d\cdot \nabla u -\tilde b\cdot \nabla v= 0,\\ \\
&\nabla\cdot v=0.
\end{array}
\right.
\end{align*}
\subsubsection{Weak solutions case}
In this section we show the uniqueness of the weak solution provided  $\alpha\geq\frac32$.  Let $(u_0,b_0), (\tilde u_0,\tilde b_0)\in L^2(\Omega) \times L^2(\Omega) $ such that  $b_0, \tilde b_0\in L^{3\frac{\alpha+1}{\alpha}}(\Omega), \nabla\cdot u_0=\nabla\cdot b_0=\nabla\cdot \tilde u_0=\nabla\cdot \tilde b_0=0, \overline b=\overline{ \tilde b}=0$. First, we multiply $\tilde{\mathcal S}_1$ by $v$ and integrate over  $\T$
\begin{align}\label{uniqv}
\nonumber\frac12\dt\vn_2^2 +\nu\vgn_2^2&+ a\,\iT \left(|u|^{2\alpha}u-\,|\tilde u|^{2\alpha}\tilde u\right) \,\cdot v\,dx \\
&= -\iT \left((v\cdot \nabla)\, \tilde u- (b\cdot \nabla)\, d-(d\cdot \nabla)\, \tilde b\right)\cdot v\,dx
\end{align}
Next, we multiply $\tilde{\mathcal S}_2$ by $d$ and integrate over  $\T$
\begin{align}\label{uniqb}
\frac12\dt\nr d\nr_2^2 +\kappa\,\nr \nabla d  \nr_2^2= \iT \left( -(v\cdot \nabla)\, b +(d\cdot \nabla)\, u +(\tilde b\cdot \nabla)\, v\right)\cdot d\,dx
\end{align}
We need the following well known monotonicity fact (see, {\it e.g.}, \cite{barrett1994finite}): There exists a positive constant $c$ such that
\begin{equation*}
c\,|u-{\tilde u}|^{2}\,\left(|u|+|{\tilde u}|\right)^{2\alpha}\leq \left(|u|^{2\alpha}u-|{\tilde u}|^{2\alpha}{\tilde u}\right)\cdot (u-{\tilde u}) .
\end{equation*}
In particular,
\begin{align}\label{positivity}
0\leq c\,a\, \iT \left(|u|+|{\tilde u}|\right)^{2\alpha} |v|^2\,dx\leq a\,\iT \left(|u|^{2\alpha}u-\,|\tilde u|^{2\alpha}\tilde u\right) \,\cdot v\,dx.
\end{align}
Now, estimate the terms appearing in the right-hand side of \eqref{uniqv} and \eqref{uniqb}. To take advantage of the velocity damping term,  we write the following
\begin{align}\label{part0uniq}
\nonumber\left|\iT (v\cdot \nabla )\,\tilde u\cdot v\,dx\right|&=\left|\iT (v\cdot \nabla )\,v\cdot \tilde u\,dx \right|\\
\text{\tiny{[By H\"older's  ineq.]}}& \leq \iT\,|\tilde u |\,|v|^\frac1\alpha\,|v|^{1-\frac1\alpha}\,|\nabla\,v|\,dx \leq \nr |\tilde u |\,|v|^\frac1\alpha\nr_{{2\alpha}}\,\nr |v|^{1-\frac1\alpha}|\nr_{\frac{2\alpha}{\alpha-1}}\,\nr \nabla v\nr_{2}\nonumber\\
\text{\tiny{[By Young's  ineq.]}}&\leq \frac 5{\nu}\, \nr|\tilde u |^{\alpha}\,|v|\nr^{\frac2\alpha}_{2}\,\nr v\nr^{2(1-\frac1\alpha)}_{2}  +\frac\nu{5}\,\nr \nabla v\nr^2_{2}\nonumber\\
\text{\tiny{[By Young's  ineq.]}}& \leq \frac a2\, \nr|\tilde u |^{\alpha}\,v\nr^{2}_{2} + \left(\frac{2\cdot 5^\alpha}{a\nu^{\alpha}}\right)^{\frac1{\alpha-1}} \,\nr v\nr^{2}_{2} +\frac\nu5\,\nr \nabla v\nr^2_{2}.
\end{align}
Now, we handle the second term in the right-hand side of \eqref{uniqv} as follows
\begin{align}\label{part1uniq}
\nonumber\left|\iT (b\cdot \nabla )\,d\cdot v\,dx\right| &=\left|\iT (b\cdot \nabla )\,v\cdot d\,dx\right| \\ \nonumber\text{\tiny{[By Cauchy-Schwarz and Young's ineq.]}}&\leq \frac \nu5\,\nr \nabla v\nr_2^2 +\frac 5\nu\,\nr b|d|\nr_2^2 \\
\nonumber\text{\tiny{[By H\"older's  ineq.]}}&\leq \frac \nu5\,\nr \nabla v\nr_2^2 +\frac 5\nu\,\nr b\nr_{3\frac{\alpha+1}{\alpha}}^2\,\nr d\nr_{6\frac{\alpha+1}{\alpha+3}}^2 \\
\text{\tiny{[By Gagliardo-Nirenberg  ineq.]}}& \leq \frac \nu5\,\nr \nabla v\nr_2^2 +\frac {c}\nu\,\nr b\nr_{3\frac{\alpha+1}{\alpha}}^2\,\nr d\nr^\frac2{\alpha+1}_2\,\nr \nabla d\nr_6^{\frac{2\alpha}{\alpha+1}}\nonumber \\
\text{\tiny{[By Young's  ineq.]}}&\leq \frac \nu5\,\nr \nabla v\nr_2^2+\frac {c}{\kappa^\alpha\nu^{\alpha+1}}\,\nr b\nr^{2({\alpha+1})}_{3\frac{\alpha+1}{\alpha}}\,\nr d\nr^2_2 + \frac\kappa6\nr \nabla d\nr_2^{2},
\end{align}
 For the last term of the right-hand side of \eqref{uniqv} we use the same argument as in the estimate above and get
\begin{flalign}\label{part2uniq}
\nonumber \left|\iT (d\cdot \nabla )\,\tilde b\cdot v\,dx\right|&=\left|\iT (d\cdot \nabla )\,v\cdot  \tilde b\,dx\right| \\\nonumber&\leq \,\frac\nu5\nr \nabla v\nr^2_2 +\frac 5\nu\nr \tilde b |d|\nr^2_2 \\
&\leq  \frac \nu5\,\nr \nabla v\nr_2^2+\frac {c}{\kappa^\alpha\nu^{\alpha+1}}\,\nr \tilde b\nr^{2({\alpha+1})}_{3\frac{\alpha+1}{\alpha}}\,\nr d\nr^2_2 + \frac\kappa6\nr \nabla d\nr_2^{2}.
\end{flalign}
The last line of the above inequality is obtained by replacing $b$ by $\tilde b$ in the estimate of $\nr bd\nr_2^2$ appearing in \eqref{part1uniq}. Next, we estimate of the terms appearing in the right-hand side of \eqref{uniqb}.
\begin{align}\label{part3uniq}
\nonumber\left|\iT (v\cdot \nabla )\,b\cdot d\,dx\right| &=\left|\iT (v\cdot \nabla d)\cdot b\,dx \right| \\
\nonumber\text{\tiny{[By Cauchy-Schwarz and Young's ineq.]}}&\leq \frac \kappa6\nr\nabla d\nr_2^2 +\frac 6\kappa\,\nr v|b|\nr_2^2 \\
\nonumber\text{\tiny{[By H\"older's  ineq.]}}&\leq \frac \kappa6\nr\nabla d\nr_2^2 +\frac {6}\kappa\,\nr b\nr_{3\frac{\alpha+1}{\alpha}}^2\,\nr v\nr_{6\frac{\alpha+1}{\alpha+3}}^2 \\
\nonumber\text{\tiny{[By Gagliardo-Nirenberg  ineq.]}}&\leq\frac \kappa6\nr\nabla d\nr_2^2 +\frac {c}\kappa\,\nr b\nr_{3\frac{\alpha+1}{\alpha}}^2\,\nr v\nr_2^{\frac2{\alpha+1}}\nr \nabla v\nr_2^{\frac{2\alpha}{\alpha+1}} + \frac {c}\kappa\,\nr b\nr_{3\frac{\alpha+1}{\alpha}}^2\,\nr v\nr_2^2 \\
\text{\tiny{[By Young's  ineq.]}}&\leq \frac\kappa6\nr\nabla d\nr_2^2 +\frac {c}{\kappa^{\alpha+1}\nu^{\alpha}}\,\nr b\nr^{2\frac{\alpha+1}{\alpha}}_{3\frac{\alpha+1}{\alpha}}\,\nr v\nr_2^{2}+ \frac\nu5 \nr \nabla v\nr_2^{2}+ \frac {c}\kappa\,\nr b\nr_{3\frac{\alpha+1}{\alpha}}^2\,\nr v\nr_2^2.
\end{align}
Next, we have
\begin{align*}
\nonumber \left|\iT (\tilde b\cdot \nabla )\,v\cdot  d\,dx\right|& \leq \left| \iT (\tilde b\cdot \nabla )\,d\cdot  v\,dx \right|\leq \frac\kappa6\,\nr \nabla d\nr_2^2 + \frac c\kappa\,\nr v|\tilde b|\nr_2^2.
\end{align*}
The above term can be estimated clearly as in \eqref{part3uniq} replacing $b$ by $\tilde b$. We get
\begin{align}\label{part6uniq}
\left| \iT (\tilde b\cdot \nabla )\,v\cdot  d\,dx\right|& \leq \frac\kappa6\nr\nabla d\nr_2^2 +\frac {c}{\kappa^{\alpha+1}\nu^{\alpha}}\,\nr \tilde b\nr^{2\frac{\alpha+1}{\alpha}}_{3\frac{\alpha+1}{\alpha}}\,\nr v\nr_2^{2}+ \frac\nu5 \nr \nabla v\nr_2^{2}+ \frac {c}\kappa\,\nr \tilde b\nr_{3\frac{\alpha+1}{\alpha}}^2\,\nr v\nr_2^2.
\end{align}
The last  term to estimate is the following
\begin {align*}
\nonumber \left|\iT (d\cdot \nabla )\,u\cdot  d\,dx\right|&=\left|\iT (d\cdot \nabla d)\cdot  u\,dx \right| \\\text{\tiny{[By Cauchy-Schwarz  ineq.]}}&\nonumber\leq  \nr \nabla d\nr_2 \,\nr u|d|\nr_{2} \\
\nonumber\text{\tiny{[By  H\"{o}lder  ineq.]}}&\leq \nr \nabla d\nr_2\,\nr u\nr_{2\alpha+2}\,\nr d\nr_{2\frac{\alpha+1}{\alpha}} \\
\nonumber\text{\tiny{[By  interpolation ]}}&\leq  \nr \nabla d\nr_2 \,\nr u\nr_{2\alpha+2}\,\nr d\nr_2^{\frac{2\alpha-1}{2\alpha+2}}\,\nr d\nr^{\frac3{2\alpha+2}}_{6}\\
\nonumber\nonumber\text{\tiny{[By  Sobolev ineq.]}}&\leq  \hat c^{\frac3{2\alpha+2}}\,\nr \nabla d\nr^{\frac{2\alpha+5}{2\alpha+2}}_2  \,\nr u\nr_{2\alpha+2}\,\nr d\nr_2^{\frac{2\alpha-1}{2\alpha+2}}\\
\text{\tiny{[By Young's  ineq.]}}&\leq  \frac \kappa3\nr \nabla d\nr^2_2+ \hat c^{\frac6{2\alpha-1}}\, \left(\frac {3}\kappa\right)^\frac{2\alpha+5}{2\alpha-1}\,\nr u\nr^{\frac{4\alpha+4}{2\alpha-1}}_{2\alpha+2}\,\nr d\nr_2^{2}.
\end{align*}
where $\hat c$ denotes the Sobolev constant. Recall that $u \in L^{2\alpha+2}_{\rm loc}([0,+\infty);L^{2\alpha+2}(\Omega))$, thanks to \eqref{propweaksolu}. Now, observe that
\[ \frac{4\alpha+4}{2\alpha-1} \leq 2\alpha+2\quad \Longleftrightarrow\quad \alpha\geq\frac32.\]
Therefore,  if $\alpha\geq \frac32$, then $u \in L^\frac{2\alpha+5}{2\alpha-1}_{\rm loc}([0,+\infty);L^{2\alpha+2}(\Omega))$  and we have
\begin{align}\label{part5uniq}
\iT (d\cdot \nabla )\,u\cdot  d\,dx&\leq  \frac \kappa3\nr \nabla d\nr^2_2+ c\, \kappa^\frac{2\alpha+5}{1-2\alpha}\,\left(1+\nr u\nr^{{2\alpha+2}}_{2\alpha+2}\right)\,\nr d\nr_2^{2}.
\end{align}
Gathering the estimates (\ref{positivity}-\ref{part5uniq})  together, we get
\begin{align*}
&\dt \left(\nr v(t)\nr_2^2+\nr d(t)\nr_2^2\right) \\ &\leq  c\,\left\{ \left({a\nu^{\alpha}}\right)^{\frac1{1-\alpha}} + \frac {1}{\kappa^{\alpha+1}\nu^{\alpha}}\,\left(\nr b\nr^{2 \frac{\alpha+1}\alpha}_{3\frac{\alpha+1}{\alpha}}+\nr \tilde b\nr^{2\frac{\alpha+1}\alpha}_{3\frac{\alpha+1}{\alpha}}\right)+  \frac {1}\kappa\,\left(\nr b\nr_{3\frac{\alpha+1}{\alpha}}^2   +\nr \tilde b\nr_{3\frac{\alpha+1}{\alpha}}^2 \right)   \right\}\,\nr v\nr^{2}_{2} \\
&+ c\,\left\{ \frac {1}{\kappa^\alpha\nu^{\alpha+1}}\,\left(\nr b\nr^{2({\alpha+1})}_{3\frac{\alpha+1}{\alpha}}+\nr \tilde b\nr^{2({\alpha+1})}_{3\frac{\alpha+1}{\alpha}}\right)+  \kappa^\frac{2\alpha+5}{1-2\alpha}\,\left(1+\nr u\nr^{{2\alpha+2}}_{2\alpha+2}\right)\right\}\,\nr d\nr^2_2.
\end{align*}
Thanks to \eqref{propweaksolu}, \eqref{bregfortimeder}, and Gronwall's inequality, we obtain clearly the continuous dependence of the weak solution on the initial data, in particular its uniqueness provided $\alpha\geq\frac32$.
\subsubsection{Strong solutions case}
In this section we show the uniqueness of the strong solution provided  $\alpha\geq\frac32$.  Let $(u_0,b_0), (\tilde u_0,\tilde b_0)\in H^1(\Omega) \times H^1(\Omega) $ such that  $\nabla\cdot u_0=\nabla\cdot b_0=\nabla\cdot \tilde u_0=\nabla\cdot \tilde b_0=0, \overline b_0=\overline{ \tilde b}_0=0$ and $\alpha\frac32$. The proof uses the estimates of the previous section except \eqref{part5uniq}. Indeed, using H\"older, Gagliardo-Nirenberg, and Young inequalities, we have
\begin {align}\label{lastforunique}
 \nonumber \iT (d\cdot \nabla u)\cdot  d\,dx \leq \nr \nabla u\nr_2 \nr d\nr_4^2 &\leq c\,\nr \nabla u\nr_2 \,\nr \nabla d\nr_2^\frac32\,\nr d\nr_2^\frac12  \\ &\leq \frac\kappa3\,\nr \nabla d\nr_2^2 + c\,{\kappa}^{-3}\,\nr \nabla u\nr_2^4 \,\nr d\nr^2_2.
\end{align}
Eventually, Summing up \eqref{uniqv} and \eqref{uniqb}, and using \eqref{positivity} along with (\ref{part0uniq}--\ref{part6uniq}) and \eqref{lastforunique}, we obtain
\begin{align*}
&\dt \left(\nr v(t)\nr_2^2+\nr d(t)\nr_2^2\right) \\ &\leq  c\,\left\{ (a\nu^{\alpha})^{\frac1{\alpha-1}} + \frac {1}{\kappa^{\alpha+1}\nu^{\alpha}}\,\left(\nr b\nr^{2 \frac{\alpha+1}\alpha}_{3\frac{\alpha+1}{\alpha}}+\nr \tilde b\nr^{2\frac{\alpha+1}\alpha}_{3\frac{\alpha+1}{\alpha}}\right)+  \frac {1}\kappa\,\left(\nr b\nr_{3\frac{\alpha+1}{\alpha}}^2   +\nr \tilde b\nr_{3\frac{\alpha+1}{\alpha}}^2 \right)   \right\}\,\nr v\nr^{2}_{2} \\
&+ c\,\left\{ \frac {1}{\kappa^\alpha\nu^{\alpha+1}}\,\left(\nr b\nr^{2({\alpha+1})}_{3\frac{\alpha+1}{\alpha}}+\nr \tilde b\nr^{2({\alpha+1})}_{3\frac{\alpha+1}{\alpha}}\right)+ \kappa^{-3}\,\nr \nabla u\nr_2^4 \right\}\,\nr d\nr^2_2.
\end{align*}
Thanks to \eqref{bregfortimeder} and \eqref{regularityub}, and Gronwall's inequality, we obtain clearly the continuous dependence of the strong solution on the initial data, in particular its uniqueness for all $\alpha\geq \frac32$.
\begin{remark}
After the completion of an earlier version of this paper we become aware of the work \cite{scooped} (see also references therein). In \cite{scooped} the authors prove   the existence and uniqueness of strong solutions to a variant of system $\mathcal S$, which contains an additional nonlinear damping term in the evolution equation of the magnetic field.  In this paper we do not need the additional nonlinear damping in the magnetic field to establish our results. Moreover, we also show  the existence and uniqueness of weak solutions. The latter has been possible thanks to the estimate that we establish in section \ref{mainsec} for the     $L^{3\frac{\alpha+1}{\alpha}}(\Omega)$ norm.
\end{remark}

\section*{Acknowledgements}
\noindent EST would like to thank the \'Ecole Polytechnique for its kind hospitality, where this work was completed, and the \'Ecole Polytechnique Foundation for its partial financial support through the 2017-2018 ``Gaspard Monge Visiting Professor'' Program.  The work of EST was also supported in part by the ONR grant N00014-15-1-2333.
\vskip6pt
\noindent ST would like to acknowledge the support provided by the Deanship of Scientific Research at King Fahd University of Petroleum \& Minerals for funding this work through project No. SR171003.

\end{document}